\pdfoutput=1 
\documentclass[10pt,twoside]{article}

\usepackage{amsmath}
\usepackage{amsfonts}
\usepackage{amssymb}
\usepackage{amscd}
\usepackage{amsthm}
\usepackage{amsbsy}
\usepackage{graphicx}
\usepackage{bm}
\usepackage{mathrsfs}

\newtheorem{Theorem}{Theorem}[section]
\newtheorem{Definition}{Definition}[section]
\newtheorem{Lemma}{Lemma}[section]

\newtheorem{Corollary}{Corollary}[section]
\newtheorem{Remark}{Remark}[section]

\def\@oddhead{\hfill \shorttitle \hfill \thepage}
\def\@evenhead{\thepage \hfill \shortauthor \hfill}
\def\@oddfoot{}
\def\@evenfoot{}

\date{}
\title{\ \\[0.4cm] \ \\ \bf  Hyperbolic Gradient Flow: Evolution of Graphs in  $\mathbb{R}^{n+1}$}
\author{De-Xing
Kong\footnote{Department of Mathematics, Zhejiang University,
Hangzhou 310027, China}\hspace{2mm} and Kefeng
Liu\footnote{Department of Mathematics, University of California at
Los Angeles, CA 90095, USA}}
\begin{document}

\maketitle


\thispagestyle{empty}

\begin{abstract}
\vskip 3mm\footnotesize{

\vskip 4.5mm \noindent In this paper we introduce a new geometric
flow --- the hyperbolic gradient flow for graphs in the
$(n+1)$-dimensional Euclidean space $\mathbb{R}^{n+1}$. This kind of
flow is new and very natural to understand the geometry of
manifolds. We particularly investigate the global existence of the
evolution of convex hypersurfaces in $\mathbb{R}^{n+1}$ and the
evolution of plane curves, and prove that, under the hyperbolic
gradient flow, they converge to the hyperplane and the straight
line, respectively, when $t$ goes to the infinity. Our results show
that the theory of shock waves of hyperbolic conservation laws can
be naturally applied to do surgery on manifolds. Some fundamental
but open problems are also given.

\vspace*{2mm} \noindent{\bf 2000 Mathematics Subject
Classification:} 53C44, 53C21, 58J45, 35L45.

\vspace*{2mm} \noindent{\bf Keywords and Phrases: }} Hyperbolic
gradient flow, geometry of manifold, global existence, smooth
solution, shock wave.

\end{abstract}

\vspace*{-12.2cm}
\vspace*{2mm} \noindent\hspace*{82mm}
\begin{picture}(41,10)(0,0)\thicklines\setlength{\unitlength}{1mm}
\put(0,2){\line(1,0){41}} \put(0,16){\line(1,0){41}}%
\put(0,12.){\sl \copyright\hspace{1mm}Higher Education Press}
\put(0,7.8){\sl \hspace*{4.8mm}and International Press}%
\put(0,3.6){{\sl \hspace*{4.8mm}Beijing--Boston} }
\end{picture}

\vspace*{-17.6mm}\noindent{{\sl The title of\\
This book*****}\\SMM\,?, pp.\,1--?} \vskip8mm
\vspace*{11.6cm}

\section{Introduction}
Classical differential geometry has been the study of curved spaces,
shapes and structures of manifolds in which the time does not play a
role. However, in the last several decades geometers have made great
strides in understanding the shapes and structures of manifolds that
evolve in time. There are many processes in the evolution of a
manifold, among them the Ricci flow is arguably the most successful
(see Hamilton \cite{h}), since it plays a fundamental role in the
solution of the famous Poincar\'{e} conjecture (see
\cite{p1}-\cite{p3}). The Ricci flow is described by a fully
nonlinear system of parabolic partial differential equations of
second order. Another famous geometric flow --- mean curvature flow
is also described by a fully nonlinear system of parabolic partial
differential equations of second order. The (inverse) mean curvature
flow has been used to prove the Riemannian-Penrose inequality in
general relativity by Huisken and Ilmanen (see \cite{hi}) and also
has been used to study many problems arising from applied fields,
i.e., imaging processing (see \cite{ak}). In fact, the traditional
geometric analysis has been successfully applied the theory of
elliptic and parabolic partial differential equations to
differential geometry and physics (see \cite{sy4}). There are three
typical examples: the Hamilton's Ricci flow, the (inverse) mean
curvature flow and the Schoen-Yau's solution of the positive mass
conjecture (see \cite{sy2}-\cite{sy3}). On the other hand, since the
hyperbolic equation or system is one of the most natural models in
the nature, a natural and important question is if we can apply the
theory of hyperbolic differential equations to solve some problems
arising from differential geometry and theoretical physics (in
particular, general relativity). Recently, we introduced the
hyperbolic geometric flow which is an attempt to answer the above
question. The hyperbolic geometric flow is a very natural tool to
understand the wave character of the metrics, wave phenomenon of the
curvatures, the evolution of manifolds and their structures (see
\cite{k}, \cite{kl}, \cite{klx}, \cite{dkl0}, \cite{dkl},
\cite{klw}, \cite{he}, \cite{kw}, \cite{kw2}).

In this paper we introduce a new geometric flow --- the hyperbolic
gradient flow for graphs in the $(n+1)$-dimensional Euclidean space
$\mathbb{R}^{n+1}$. The flow is described by hyperbolic evolution
partial differential equations of first order for a family of vector
fields $X_t$ defined on $\mathbb{R}^{n}$. Roughly speaking, the
hyperbolic gradient flow evolves the tangent planes of the graph
under consideration, this is different from the Ricci flow, the mean
curvature flow or our hyperbolic geometric flow. This kind of flow
is new and very natural to understand deformation phenomena of
manifolds (in particular, graphs in $\mathbb{R}^{n+1}$) as well as
the geometry of manifolds. It possesses many interesting properties
from both mathematics and physics. In the present paper, we
particularly investigate the global existence of the evolution of
convex hypersurfaces in $\mathbb{R}^{n+1}$ and the evolution of
plane curves, and prove that, under the hyperbolic gradient flow,
they converge to the hyperplane and the straight line, respectively,
when $t$ goes to the infinity. Our results show that the theory of
shock waves of hyperbolic conservation laws can be naturally applied
to do surgery on manifolds. Some fundamental but open problems are
also given.


\section{Hyperbolic gradient flow for graphs in $\mathbb{R}^{n+1}$}

Let $\Sigma_t$ be a family of graphs in the $(n+1)$-dimensional
Euclidean space $\mathbb{R}^{n+1}$ with coordinates $(x_1,\cdots,
x_{n+1})$. Without loss of generality, we may assume that the graphs
$\Sigma_t$ are given by
\begin{equation}x_{n+1}=f(t,x_1,\cdots,x_n),\end{equation}
where $f$ is a smooth function defined on
$\mathbb{R}\times\mathbb{R}^{n}$. Let $X_t$ be a family of tangent
vector fields induced by $\Sigma_t$, or say,
\begin{equation} X_t=(X_1,\cdots X_n)=(\partial_{x_1}f,\cdots, \partial_{x_n}f),\end{equation}
where $\partial_{x_i}f\;(i=1,\cdots,n)$ stand for $\frac{\partial
f}{\partial x_i}$. The hyperbolic gradient flow under considered
here is given by the following evolution equations
\begin{equation}\label{2.1}\frac{\partial X_t}{\partial t}+\nabla\left(\frac{\|X_t\|^2}{2}\right)=0,\end{equation}
where $\nabla =(\partial_{x_1},\cdots,
\partial_{x_n})$ and
\begin{equation}\|\cdot\|^2=\langle\cdot,\cdot\rangle,\end{equation}
in which $\langle\cdot,\cdot\rangle$ stands for the inner product in
$\mathbb{R}^{n}$.

By the definition, the hyperbolic gradient flow introduced in this
note is a geometric flow for the evolution of a family of tangent
vector fields induced by a family of graphs, it is quite different
from the Ricci flow and the mean curvature flow: the Ricci flow is
described by evolution equations for a family of Riemannian metrics
$g_{ij}(t)$ defined on the manifold under consideration, while the
mean curvature flow is on the evolution of the manifold itself.

\section{The evolution of convex hypersurfaces in  $\mathbb{R}^{n+1}$}
In this section, we shall investigate the evolution of convex
hypersurfaces in the $(n+1)$-dimensional Euclidean space
$\mathbb{R}^{n+1}$.

As before, let $\mathbb{R}^{n+1}$ be the $(n+1)$-dimensional
Euclidean space with coordinates $(x_1,\cdots, x_{n+1})$, and
$x_{n+1}=\mathscr{S}(t,x_1,\cdots, x_n)$ be a family of
hypersurfaces in $\mathbb{R}^{n+1}$. Introduce the vector field
\begin{equation}\label{3.1}
\vec{v}=\{\mathscr{S}_1,\cdots,\mathscr{S}_n\},\end{equation} where
$\mathscr{S}_i\;(i=1,\cdots,n)$ stand for $\frac{\partial
\mathscr{S}}{\partial x_i}$. In the present situation, the
hyperbolic gradient flow (\ref{2.1}) reads
\begin{equation}\label{3.2}
\vec{v}_t+\vec{v}\cdot\nabla\vec{v}=0.\end{equation} In this case,
(\ref{3.2}) is nothing but the transport equation for $\vec{v}$.

\vskip 4mm

\noindent {\bf Example 3.1.} {\it  Consider the evolution of the
hypersurface $x_n=\frac12\left(x_1^2+\cdots+x_n^2 \right)$ under the
hyperbolic gradient flow. In the present situation, we need to
consider the Cauchy problem for the equation (\ref{3.2}) with the
following initial data
\begin{equation}\label{3.3}
t=0:\;\; \vec{v}=\vec{v}^0\triangleq (x_1,\cdots,x_n).\end{equation}
It is easy to see that the solution of the Cauchy problem
(\ref{3.2}), (\ref{3.3}) reads
\begin{equation}\label{3.4}
\vec{v}=\left(\frac{x_1}{t+1},\cdots,\frac{x_n}{t+1}\right),\end{equation}
moreover, the solution is unique. Obviously, the vector field
defined by (\ref{3.4}) gives a potential function
$x_n=\frac{1}{2(t+1)}\left(x_1^2+\cdots+x_n^2\right)+C$, where $C$
is a constant independent of $x$. Noting that the initial
hypersurface is $x_n=\frac12\left(x_1^2+\cdots+x_n^2 \right)$ leads
to that the constant $C$ must be zero. Thus, the evolution of the
hypersurface $x_n=\frac12\left(x_1^2+\cdots+x_n^2 \right)$ under the
hyperbolic gradient flow is described by the family of hypersurfaces
$x_n=\frac{1}{2(t+1)}\left(x_1^2+\cdots+x_n^2\right)$. Clearly, for
any fixed $x$, the hypersurfaces tend to flat under the hyperbolic
gradient flow when $t$ goes to the infinity.} $\qquad\qquad\Box$

Consider the Cauchy problem for the equation (\ref{3.2}) with the
following initial data
\begin{equation}\label{3.5}
t=0:\;\; \vec{v}=\vec{v}^0(x_1,\cdots,x_n),\end{equation} where
$\vec{v}^0$ is a smooth vector field defined on $\mathbb{R}^{n}$. We
now consider the global existence and decay property of smooth
solutions of the the Cauchy problem (\ref{3.2}) and (\ref{3.5}).

In fact, we can obtain a sufficient and necessary condition on the
global existence of smooth solutions of the following Cauchy problem
for more general quasilinear systems of first order
\begin{equation}\label{1.3.1}
\frac{\partial u}{\partial t}+\sum^n_j \lambda_j(u)\frac{\partial
u}{\partial x_j}=0,\quad \forall\; (t, x)\in \mathbb{R}^+\times
\mathbb{R}^n
\end{equation}
with the initial data
\begin{equation}\label{1.3.2}
u(0,x)=\phi(x), \quad \forall\; x\in \mathbb{R}^n,
\end{equation}
where $x=(x_1,\cdots ,x_n)$ stands for the special variable,
$u=(u_1(x, t),\cdots  u_m(x,t))^T$ is the unknown vector-valued
function of $(t, x)=(t, x_1,\cdots  ,x_n)\in \mathbb{R}^+\times
\mathbb{R}$, $\lambda_i(u)\; (i=1,\cdots , n)$ are given $C^1$
functions and $\phi(x)=(\phi_1(x),\cdots , \phi_m(x))^T$ is a given
$C^1$ vector-valued function with bounded $C^1$ norm. The following
lemma comes from Conway \cite{conway}, Li \cite{li}, Dafermos
\cite{dafermos} or Kong \cite{kong}.

\begin{Lemma} Under the assumptions mentioned above, the
Cauchy problem (\ref{1.3.1})-(\ref{1.3.2}) has a unique global $C^1$
smooth solution on the domain $\mathbb{R}^+\times \mathbb{R}^n$ if
only if, for any given $x\in \mathbb{R}^n$, it holds that
\begin{equation}\label{1.3.3}
d(S_p V_0(x), \mathbb{R}^- )\geq0,
\end{equation}
i.e., all eigenvalues of the $n\times n$ matrix
\begin{equation}\label{1.3.4}
V_0(x)=\left(\sum^m_{k=1}\frac{\partial \lambda_i}{\partial
u_k}(\phi(x))\frac{\partial \phi_k}{\partial x_j}\right)^n_{i,j=1}
\end{equation}
are non-negative, where $S_p V_0(x)$ stands for the spectrum of the
matrix $V_0(x)$.
\end{Lemma}

\begin{Lemma}
Under the assumptions of Lemma 3.2, suppose that $\phi$ is a $C^2$
vector-valued function with bounded $C^2$ norm and suppose
furthermore that there exists a positive constant $\delta > 0$ such
that
\begin{equation}\label{1.3.5}
d(S_p V_0(x), \mathbb{R}^- )\geq \delta, \quad \forall\; x\in
\mathbb{R}^n.
\end{equation}
Then the Cauchy problem (\ref{1.3.1})-(\ref{1.3.2}) admits a unique
global $C^2$ smooth solution $u=u(t,x)$ on the domain
$\mathbb{R}^+\times \mathbb{R}^n$, moreover it holds that
\begin{equation}\label{1.3.6}
\parallel Du(t,x)\parallel_{L^{\infty}(\mathbb{R}^n)}=C_1(1+t)^{-1}
\end{equation}
and
\begin{equation}\label{1.3.7}
\parallel D^2u(t,x)\parallel_{L^{\infty}(\mathbb{R}^n)}\leq
C_2(1+t)^{-2},
\end{equation}
where $C_1$ is a positive constant independent of $t$ but depending
on $\delta$ and the $C^1$ norm of $\phi$, while $C_2$ is a positive
constant independent of $t$ but depending on $\delta$ and the $C^2$
norm of $\phi$.
\end{Lemma}

The proof of Lemma 3.3 can be found in  Grassin \cite{g} for the
case of scalar equation and in Kong \cite{kong} for general case.

\vskip 4mm

If $m=n$ and $\lambda_i(u)=u_i$ $(i=1,\cdots,n)$ (in this case, the
system (\ref{1.3.1}) goes back to the system (\ref{3.2})), then in
the present situation, $V_0(x)$ defined by (\ref{1.3.4}) reads
\begin{equation}\label{1.3.38}
V_0(x)=\left(\frac{\partial \phi_i}{\partial x_j}\right)^n_{i,j=1}.
\end{equation}
In particular, if there exists a potential function $\Phi(x)$ such
that
\begin{equation}\label{1.3.39}
\frac{\partial \Phi}{\partial x_i}=\phi_i(x) \quad (i=1,\cdots, n),
\end{equation}
then
\begin{equation}\label{1.3.40}
V_0(x)={\rm Hess}\,(\Phi(x)).
\end{equation}

We now turn to consider the Cauchy problem for this special case,
i.e.,
\begin{equation}\label{1.3.41} \left\{\begin{array}{l}{\displaystyle
\frac{\partial u}{\partial t}+\sum^n_{j=1}u_j\frac{\partial
u}{\partial x_j}=0,\quad
\forall\; (t,x)\in \mathbb{R}^+\times \mathbb{R}^n,}\vspace{2mm}\\
{\displaystyle t=0: u=\phi(x)=\left(\frac{\partial \Phi}{\partial
x_1},\cdots , \frac{\partial \Phi}{\partial x_n}\right)^T, \quad
\forall\; x\in \mathbb{R}^n.}
 \end{array}\right.\end{equation}
By Lemma 3.1-3.2, we have
\begin{Lemma} Suppose that the potential function $\Phi=\Phi(x)$ is
$C^2$ smooth and its derivates $\Phi_{k_i}$ $(i=1,..., n)$ has a
bounded $C^1$ norm. Then the Cauchy problem  (\ref{1.3.41}) has a
unique global $C^1$ smooth solution on the domain
$\mathbb{R}^+\times \mathbb{R}^n$ if only if ${\rm Hess}\,(\Phi(x))$
is non-negative for all $x\in \mathbb{R}^n$.

Moreover, if the following assumptions are satisfied: (i) $\Phi$ is
a $C^3$ smooth function, (ii) the derivative
$D\Phi=(\Phi_{x_1},\cdots, \Phi_{x_n})^T$ is a vector-valued
function with bounded $C^2$ norm, (iii) there exists a positive
constant $\delta$ independent of $x$ such that
\begin{equation}\label{1.3.42}
d({\rm Hess}\,(\Phi(x),  \mathbb{R}^-)\geq \delta, \quad  \forall\;
x\in \mathbb{R}^n,
\end{equation}
then the global smooth solution $u=u(t,x)$ to the Cauchy problem
(\ref{1.3.41}) satisfies the following properties:

(I) there exists a $C^3$ potential function $U=U(t,x)$ such that
\begin{equation}\label{1.3.43}
u_i(t,x)=\frac{\partial U}{\partial x_i}(t,x)  \quad (i=1, ..., n),
\quad \forall\; (t,x)\in \mathbb{R}^+\times \mathbb{R}^n,
\end{equation}

(II) there exist two positive constants $C_3$ and $C_4$ independent
of $t$ but depending on $\delta$ and the $C^1$ norm (for $C_3$), the
$C^2$ norm (for  $C_4$) of $D\Phi(x)$, respectively, such that
\begin{equation}\label{1.3.44}
\| D^2 U(t, \cdot)\|_{L^{\infty}(\mathbb{R}^n)}\leq C_3(1+t)^{-1},
\end{equation}
and
\begin{equation}\label{1.3.45}
\| D^3 U(t, \cdot)\|_{L^{\infty}(\mathbb{R}^n)}\leq C_4(1+t)^{-2},
\end{equation}
where $D=(\partial_{x_1}, ..., \partial_{x_n})$.
\end{Lemma}
\noindent{\bf Proof.}
 By Lemmas 3.1-3.2, we only need to prove (I) in Lemma 3.3.

In order to prove (I), it suffices to show
\begin{equation}\label{1.3.46}
\frac{\partial u_i(t,x)}{\partial x_j}=\frac{\partial
u_j(t,x)}{\partial x_i}, \quad \forall\; i\neq j,  \quad \forall\;
(t,x) \in \mathbb{R}^+\times \mathbb{R}^n.
\end{equation}

In fact, introduce
\begin{equation}\label{1.3.47}
\omega^i_j=\frac{\partial u_i}{\partial x_j}-\frac{\partial
u_j}{\partial x_i}  \quad (i,j=1,\cdots,n;\quad i \neq j) .
\end{equation}
Obviously, when $t=0$,
\begin{equation}\label{1.3.48}
\omega^i_j(t,0)=\frac{\partial}{\partial
x_j}\left(\frac{\partial}{\partial
x_i}\Phi\right)-\frac{\partial}{\partial
x_i}\left(\frac{\partial}{\partial x_j}\Phi\right)=0 \quad
(i,j=1,\cdots,n).
\end{equation}
On the one hand, differentiating the $i$-th equation in
(\ref{1.3.41}) with respect to $x_j$ gives
\begin{equation}\label{1.3.49}
\frac{\partial}{\partial t}\left(\frac{\partial u_i}{\partial
x_j}\right)+ \sum^n_{k=1} u_k \frac{\partial}{\partial
x_k}\left(\frac{\partial u_i}{\partial x_j}\right)= -\sum^n_{k=1}
\frac{\partial u_i}{\partial x_k}\frac{\partial u_k}{\partial x_j}.
\end{equation}
On the other hand, differentiating the $j$-th equation in
(\ref{1.3.41}) with respective to $x_i$ yields
\begin{equation}\label{1.3.50}
\frac{\partial}{\partial t}\left(\frac{\partial u_j}{\partial
x_i}\right)+ \sum^n_{k=1} u_k \frac{\partial}{\partial
x_k}\left(\frac{\partial u_j}{\partial x_i}\right)=
-\sum^n_{k=1}\frac{\partial u_j}{\partial x_k}\frac{\partial
u_k}{\partial x_i}.
\end{equation}
Combing (\ref{1.3.49})-(\ref{1.3.50}) leads to
\begin{equation}\label{1.3.51}
\frac{\partial \omega^i_j}{\partial t}+ \sum^n_{k=1}
u_k\frac{\partial \omega^i_j}{\partial x_k}=\sum_{p \neq q}
\Gamma^{ij}_{pq} \omega ^q_p, \quad  \forall\; i \neq j,
\end{equation}
where $\Gamma^{ij}_{pq}$ stands for the coefficients of $\omega
^q_p$ which are smooth functions of $\frac{\partial u_l}{\partial
x_h}~(l,h=1 ,\cdots, n)$. Clearly, $\omega^i_j=0 ~(i,j=1,\cdots,n;
~i\neq j)$ is a solution of the Cauchy problem (\ref{1.3.51}),
(\ref{1.3.48}). By the uniqueness of the smooth solution of the
Cauchy problem for hyperbolic partial differential equation, we have
\begin{equation}\label{1.3.52}
\omega ^i_j \equiv 0 \quad (i\neq j),\quad \forall\; (t,x) \in
\mathbb{R}^+\times \mathbb{R}^n.
\end{equation}
This proves (\ref{1.3.46}). Thus the proof of Lemma 3.3 is
completed.$\qquad\qquad \blacksquare$

\begin{Remark}In Lemma 3.2, we need the $C^1$ norm of $\phi$ and the $C^2$
norm of $\phi$ is bounded for the estimates (\ref{1.3.6}) and
(\ref{1.3.7}), respectively. For Lemma 3.3, the situation is
similar. \end{Remark}

However, in many cases (i.g., Example 3.1), the assumption that the
$C^1$ norm or $C^2$ norm of the initial data is bounded is not
satisfied. The following discussion is devoted to the case of
unbounded initial data. For simplicity, we only consider the Cauchy
problem (\ref{1.3.41}).
\begin{Lemma} Suppose that
$\Phi=\Phi(x)$ is a $C^3$ convex function, i.e., $\Phi(x) \in C^3
(\mathbb{R}^n)$ and
\begin{equation}\label{1.3.53}
{\rm Hess}\,(\Phi)\geq 0.
\end{equation}
Then the Cauchy problem (\ref{1.3.41}) admits a unique $C^2$
solution $u=u(t,x)$ on the domain $\mathbb{R}^+\times \mathbb{R}^n$.
Moreover, there exists a potential function $U=U(t,x)\in C^3(
\mathbb{R}^+\times \mathbb{R}^n)$ such that (\ref{1.3.43}) is
satisfied. In particular, if there exists a positive constant
$\delta$ independent of $x$ such that (\ref{1.3.42}) holds, then for
any fixed $\alpha \in \mathbb{R}^n$ along the characteristic curve
$x=x(t,\alpha)$ it holds that
\begin{equation}\label{1.3.54}
|D^2U(t, x(t, \alpha))|\leq \tilde{C_1}(1+t)^{-1}
\end{equation}
 and
\begin{equation}\label{1.3.55}
|D^3U(t, x(t, \alpha))|\leq \tilde{C_2}(1+t)^{-2},
\end{equation}
 where $\tilde{C_1}$ and $\tilde{C_2}$ are tow constants independent of $t$ but depending
 on $\delta$ and $\alpha$.
\end{Lemma}

The following corollary comes from Lemma 3.4 directly.
\begin{Corollary} Under the assumptions of Lemma 3.4, for any
compact set $\Omega \subseteq \mathbb{R}^n$ it holds that
\begin{equation}\label{1.3.56}
\|D^2U(t, \cdot)\|_{L^{\infty}(\Omega(t))}\leq
\tilde{C_3}(1+t)^{-1},
\end{equation}
and
\begin{equation}\label{1.3.57}
\|D^3U(t, \cdot)\|_{L^{\infty}(\Omega(t))}\leq
\tilde{C_4}(1+t)^{-2},
\end{equation}
where
\begin{equation}\label{1.3.58}
\Omega(t)=\{(t,x)| x=x(t, \alpha),\quad \alpha \in \Omega \},
\end{equation}
$\tilde{C_3}$ and $\tilde{C_4}$ are two constants independent of $t$
but depending on $\delta$ and the set $\Omega$. $\qquad \Box$
\end{Corollary}

\noindent{\bf Proof of Lemma 3.4.} Noting (\ref{1.3.53}), we have
\begin{equation}\label{1.3.59}
\Phi_{x_ix_i}(x)\geq 0 \quad (i=1,\cdots, n), \quad \forall\; x \in
 \mathbb{R}^n.
\end{equation}
In the present situation, the characteristic curve pasting through
any fixed point $(0,\alpha)$ in the initial hyperplane $t=0$ reads
\begin{equation}\label{1.3.60}
x_i=\alpha_i + \frac{\partial \Phi}{\partial \alpha_i}(\alpha)t
\quad (i=1,\cdots, n).
\end{equation}
By (\ref{1.3.59}), it is easy to check that the mapping $\Pi_t:
\mathbb{R}^n \rightarrow  \mathbb{R}^n $ defined by (\ref{1.3.60})
is {\it proper}. On the other hand,
\begin{equation}\label{1.3.61}
J(\Pi_t)=I+t{\rm Hess}\,(\Phi).
\end{equation}
Using (\ref{1.3.53}) again, we have
\begin{equation}\label{1.3.62}
\det J(\Pi_t) \geq 1, \quad \forall\; (t,\alpha) \in
\mathbb{R}^+\times \mathbb{R}^n.
\end{equation}
This implies that for any fixed $x \in \mathbb{R}^+$, the mapping
$\Pi_t$ is a global $C^1$ deffeomorphism. Solving $\alpha$ from
(\ref{1.3.60}) gives
\begin{equation}\label{1.3.63}
\alpha=\alpha(t,x) \in C^2( \mathbb{R}^+\times \mathbb{R}^n).
\end{equation}
The rest of the proof is standard (See \cite{li} or \cite{kong}),
here we omit the details. The proof of Lemma 3.4 is
completed.$\qquad\qquad \blacksquare$

\vskip 4mm

Lemma 3.3 guarantees that, if the initial vector field is induced by
a graph, then so does the solution vector-field. That is, if there
exists a fuction $\varphi_0(x_1,\cdots,x_n)$ such that
$v^0_i=\frac{\partial\varphi_0}{
\partial x_i}\;(i=1,\cdots,n)$, then there is a family of functions
$\varphi(t,x_1,\cdots,x_n)$ such that
\begin{equation}\label{3.16}
v_i(t,x_1\cdots,x_n)=\frac{\partial\varphi}{\partial
x_i}(t,x_1\cdots,x_n)\quad (i=1,\cdots,n)
\end{equation}
and
\begin{equation}\label{3.17}
\varphi(0,x_1,\cdots,x_n)=\varphi_0(x_1,\cdots,x_n).
\end{equation}
From the point of view of geometry, the hyperbolic gradient flow
evolves a graph as a family of graphs in the Euclidean space
$\mathbb{R}^{n+1}$.

\vskip 4mm

Summarizing the above argument leads to
\begin{Theorem}
For any given initial vector field induced by a convex graph
$x_{n+1}=\varphi_0(x_1,\cdots,x_n)$, the solution
$\vec{v}=\vec{v}(t,x_1,\cdots,x_n)$ to the hyperbolic gradient flow
(\ref{3.2}) exists for all time, and there exists a unique family of
graphs $x_{n+1}=\varphi(t,x_1,\cdots,x_n)$ such that the solution
vector-field $\vec{v}=\vec{v}(t,x_1,\cdots,x_n)$ is induced by the
family of graphs $x_{n+1}=\varphi(t,x_1,\cdots,x_n)$. Moreover, if
the initial graph is strictly convex, then for any fixed point
$(x_1,\cdots,x_n)\in \mathbb{R}^{n}$ the graphs
$x_{n+1}=\varphi(t,x_1,\cdots,x_n)$ tends to be flat at an algebraic
rate $(t+1)^{-1}$, when $t$ goes to the infinity.
\end{Theorem}

\section{The evolution of plane curves}
In this section, we particularly investigate the evolution of plane
curves under the hyperbolic gradient flow, here we still consider
the graph case, however we do not assume that the graph is convex.

Let $y=f(x)$ be a smooth curve in the $(x,y)$-plane, and
\begin{equation}\label{4.1}
v_0(x)=f^{\prime}(x)\end{equation} be the slope function of the
curve. In the present situation, the hyperbolic gradient flow
equation (\ref{3.2}) becomes one-dimensional case, i.e.,
\begin{equation}\label{4.2}
v_t+vv_x=0.\end{equation} This equation can be rewritten as a
conservative form
\begin{equation}\label{4.3}
v_t+(v^2/2)_x=0.\end{equation}

We next consider the Cauchy problem for the conservation law
(\ref{4.3}) with the initial data
\begin{equation}\label{4.4}
t=0:\;\; v=v_0(x).\end{equation}

As in Lax \cite{lax}, we introduce
\begin{Definition}
A function $\psi$ has mean value $M$, if
\begin{equation}\label{4.5}
\lim_{L\rightarrow\infty}\frac1L\int_{a}^{a+L}\psi(x)dx=M
\end{equation}
uniformly in $a$.
\end{Definition}

\begin{Corollary}
If a function $\psi$ is periodic with $p$ period, then it has mean
value $M$, and $M$ is given by
\begin{equation}\label{4.5-1}
M=\frac1p\int_{0}^{p}\psi(x)dx.\end{equation} If $\psi\in
L^1(\mathbb{R})$, then it has mean value 0.
\end{Corollary}

The following lemma comes from Lax \cite{lax}.
\begin{Lemma}
Let $v(t,x)$ be a bounded weak solution of the Cauchy problem
(\ref{4.3}), (\ref{4.4}). Suppose that the initial data $v_0(x)$ has
a mean value, then $v(t,x)$ has the same mean value for all $t$.
\end{Lemma}

The following important lemma comes from Kruzkov \cite{kr}.
\begin{Lemma}
Suppose that the initial data $v_0$ is bounded measurable, then the
Cauchy problem (\ref{4.3}), (\ref{4.4}) has a unique entropy
solution $v=v(t,x)$ on the half plane $t\ge 0$.
\end{Lemma}

\begin{Lemma} Under the assumption of Lemma 4.2, if the initial data $v_0$
is periodic with $p$ period, then the entropy solution $v=v(t,x)$ of
the Cauchy problem (\ref{4.3}), (\ref{4.4}) tends to $M$ uniformly
in $x$ at an algebraic rate $(t+1)^{-1}$, when $t$ tends to
infinity, where $M$ is given by
$$M=\frac1p\int_{0}^{p}v_0(x)dx.$$
\end{Lemma}

\begin{Lemma}
Suppose that the initial data $v_0$ is in the class of
$L^{1}(\mathbb{R})$, then the Cauchy problem (\ref{4.3}),
(\ref{4.4}) has a unique entropy solution $v=v(t,x)$ on the half
plane $t\ge 0$. Moreover, $v(t,x)$ tends to $0$ uniformly in $x$ at
an algebraic rate $(t+1)^{-1}$, when $t$ tends to infinity.
\end{Lemma}

Lemmas 4.3 and 4.4 can be found in Serre \cite{serre} and Bressan
\cite{bressan}, respectively.

\begin{Remark} The entropy solution $v=v(t,x)$ mentioned in Lemmas 4.2, 4.3 and 4.4 means
that (i) $v=v(t,x)$ is a weak solution of the Cauchy problem
(\ref{4.3}), (\ref{4.4}); (ii) it satisfies the entropy condition.
In fact, the entropy solution may includes shock waves, rarefaction
waves, and other physical discontinuities.
\end{Remark}

We now consider the evolution of the initial curve $y=f(x)$ under
the hyperbolic gradient flow.

Without loss of generality, we may assume that
\begin{equation}\label{4.6}
f(0)=0.
\end{equation}
By Lemmas 4.2-4.4, we have
\begin{Theorem}
Suppose that $f^{\prime}(x)$ is bounded measurable, and suppose
furthermore that $f^{\prime}(x)$ is periodic or is in the class of
$L^{1}(\mathbb{R})$, then the family of curves $y=F(t,x)$ tends to
the straight line $y=Mx$ uniformly in $x$ at an algebraic rate
$(t+1)^{-1}$, when $t$ tends to infinity, where $M$ is the mean
value of $f^{\prime}(x)$, and $y=F(t,x)$ is generated by the
hyperbolic gradient flow, i.e., $F(t,x)$ satisfies
\begin{equation}\label{4.6}
\frac{\partial F}{\partial x}(t,x)=v(t,x),
\end{equation}
in which $v=v(t,x)$ is the entropy solution of the Cauchy problem
(\ref{4.3}), (\ref{4.4}) (in the present situation,
$v_0(x)=f^{\prime}(x)$).
\end{Theorem}

\begin{Remark} In geometry, one is, in general, interested in
the case that the initial data $f(x)$, or say $v_0$, is smooth and
bounded. However, in the evolution process under the hyperbolic
gradient flow, discontinuities may appear. See Example 4.1 below for
the details.
\end{Remark}

\noindent {\bf Example 4.1.} {\it Consider the evolution of the
curve $y=-\cos x$ in the $(x,y)$-plane under the hyperbolic gradient
flow. In the present situation, the initial data reads
\begin{equation}\label{4.7}
t=0:\;\; v=v_0(x)=\sin x.\end{equation} By the method of
characteristics, the solution of the Cauchy problem (\ref{4.3}),
(\ref{4.7}) can be constructed and is given by Figure 4.1.
\begin{figure}[!ht]
    \begin{center}
\vspace{6.0cm}
\hspace{-6.0cm}\scalebox{0.8}[0.55]{\includegraphics[trim=10mm 5mm
-30mm 0mm,  width=2.5cm, height=5cm,bb=10 20 170 340]{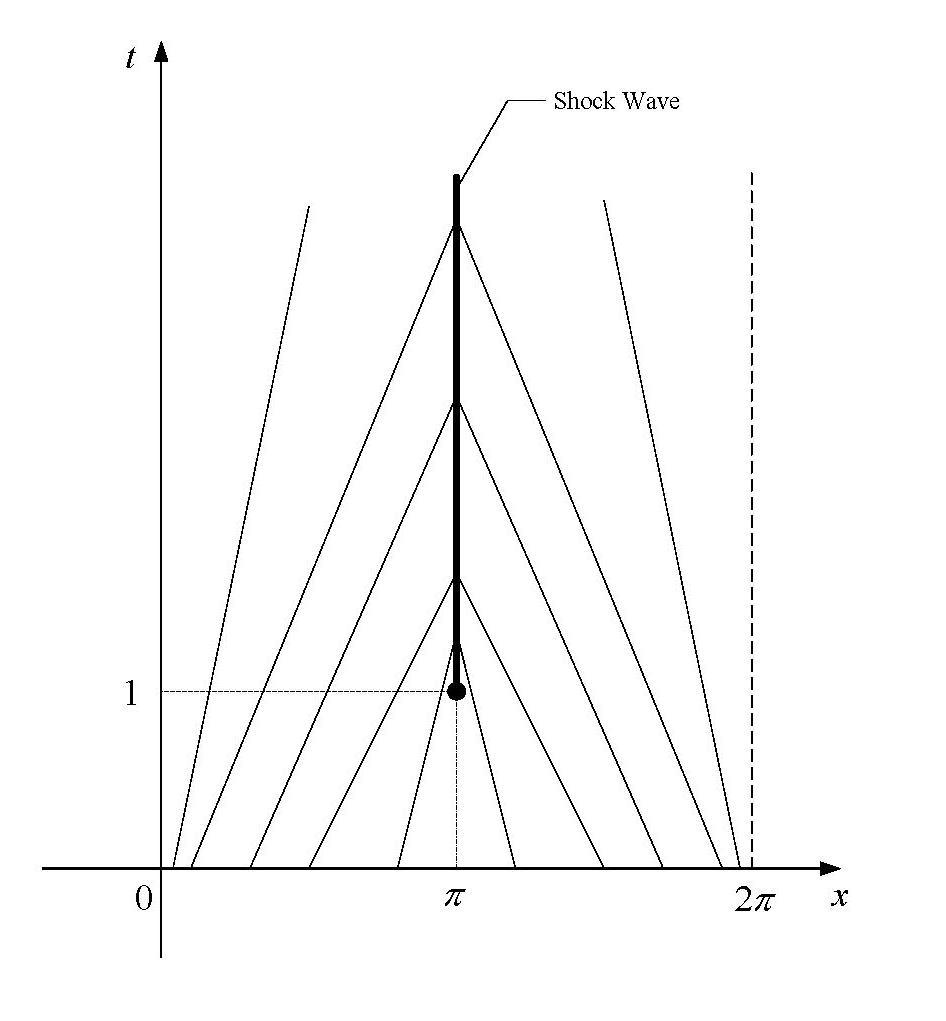}} \vskip
-2.0cm \center{Figure 4.1:~The shock solution $v=v(t,x)$ of the
Cauchy problem (\ref{4.3}), (\ref{4.7})}
    \end{center}
\end{figure}

\noindent Notice that Figure 4.1 only describes the solution on one
space-periodic domain, i.e., $ \mathbb{R}^+\times [0,2\pi]$.
Corresponding to the solution shown in Figure 4.1, the evolution of
the curve $y=-\cos x$ under the hyperbolic gradient flow can be
described by Figure 4.2.
\begin{figure}[htbp]
    \begin{center}
\vspace{7.0cm}
\hspace{-13.0cm}\scalebox{0.8}[0.55]{\includegraphics[trim=10mm 5mm
-30mm 0mm,  width=1.8cm, height=3cm,bb=10 20 170 340]{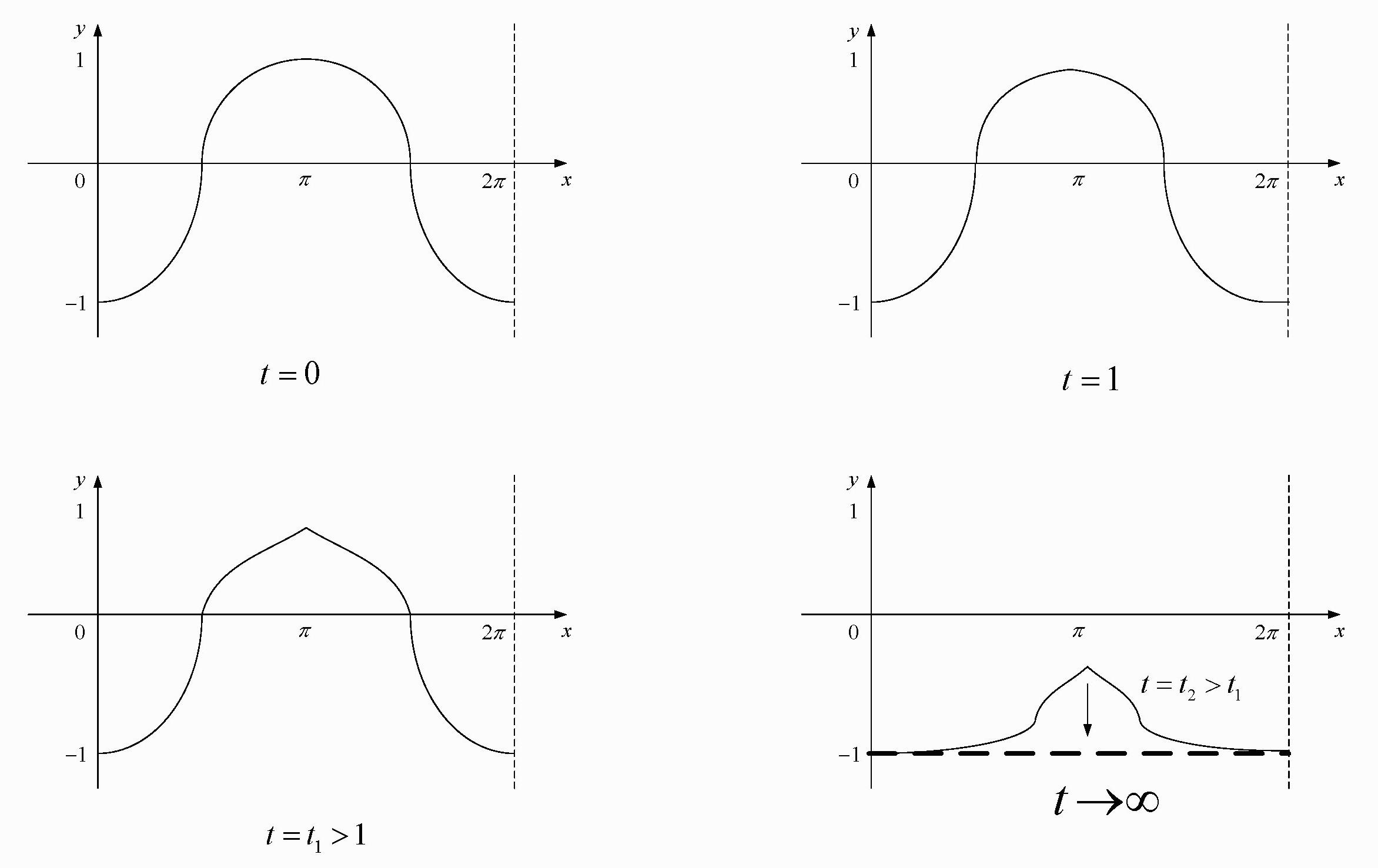}} \vskip
-2.0cm \center{Figure 4.2:~ The evolution of the curve $y=-\cos x$
under the hyperbolic gradient flow}
    \end{center}
\end{figure}

\noindent We observe from Figure 4.2 that the singularity have
appeared in the evolutionary process (see Figure 4.3 for the
details).
\begin{figure}[!h]
    \begin{center}
\vspace{3.0cm}
\hspace{-7.0cm}\scalebox{0.8}[0.55]{\includegraphics[trim=10mm 5mm
-30mm 0mm,  width=2.8cm, height=5cm,bb=10 20 170 340]{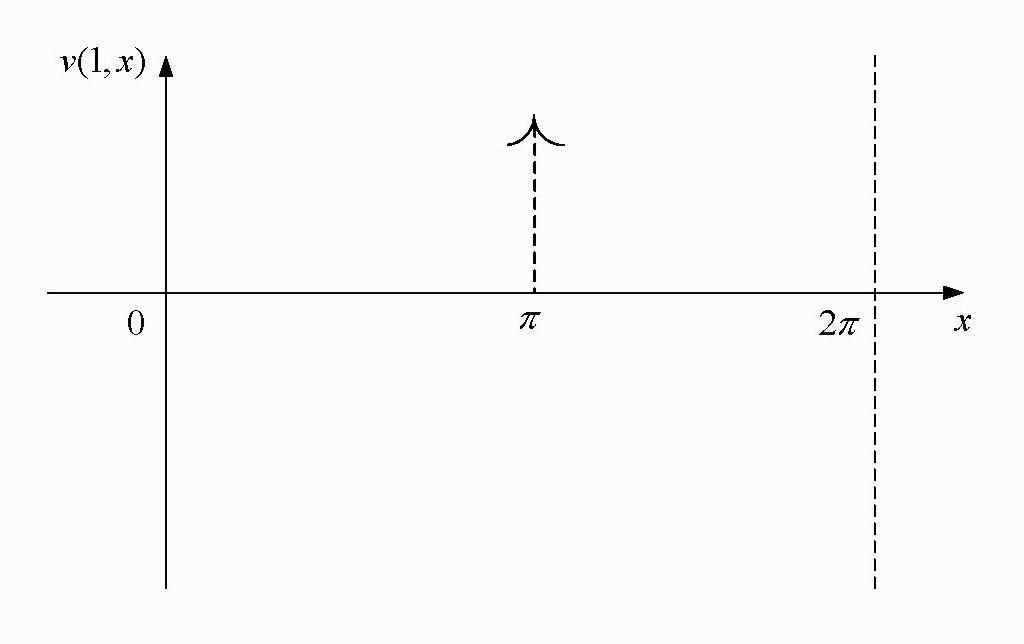}} \vskip
-2.0cm \center{Figure 4.3:~ The formation of singularity of cusp
type of $v=v(t,x)$}
    \end{center}
\end{figure}

\noindent Figure 4.3 shows that the singularity of cusp type of
$v=v(t,x)$ appears at $x=\pi$ when $t=1$. It is easy to see that the
entropy solution $v=v(t,x)$ to the Cauchy problem (\ref{4.3}),
(\ref{4.7}) includes space-periodic shock waves.} $\qquad\qquad\Box$

\section{Conclusions and open problems}

It is well known that there have been many successes of elliptic and
parabolic equations applied to mathematics and physics. On the other
hand, hyperbolic partial differential equation is a very important
kind of PDEs, it can be used to describe the wave phenomena in the
nature and engineering. Recently, we introduced the hyperbolic
geometric flow, showed that the hyperbolic geometric flow possesses
very interesting geometric properties and dynamical behavior, and
obtain some interesting results. However, the hyperbolic geometric
flow is described by a fully nonlinear system of hyperbolic partial
differential equations of second order, which is very difficult to
solve. In this paper we introduce a new geometric flow --- the
hyperbolic gradient flow, which is described by a quasilinear system
of hyperbolic partial differential equations of first order.
Comparing the hyperbolic geometric flow, the hyperbolic gradient
flow is easier to solve. The key point of the hyperbolic gradient
flow is to evolve the tangent planes of the graphs under
consideration, this is different with the famous Ricci flow, the
mean curvature flow or our hyperbolic geometric flow. In this paper,
we investigate the evolution of convex hypersurfaces in the
$(n+1)$-dimensional Euclidean space $\mathbb{R}^{n+1}$ and the
evolution of plane curves, and prove that, under the hyperbolic
gradient flow, they converge to the hyperplane and straight line,
respectively, when $t$ goes to the infinity. Our results obtained in
this paper show that the theory of shock waves of hyperbolic
conservation laws can be naturally applied to differential geometry.
We believe that the hyperbolic gradient flow is a new and powerful
tool to study some problems arising from geometry and physics.
However, there are many fundamental but still open problems. In
particular, the following open problems seem to us more interesting
and important:

{\bf 1. The evolution of plane curves.} {\it In Theorem 4.1 if we do
not assume that $f^{\prime}(x)$ has a mean value, what is the limit
of the family of curves $F(t,x)$ as $t$ goes to infinity? Moreover,
what happens if the initial curve is not a graph, e.g., a closed
curve?}

{\bf 2. The evolution of surfaces in $\mathbb{R}^3$.} {\it In
Theorem 3.1 if the initial surface is a graph but is not convex,
what about the limit of the family of surfaces $\varphi(t,x_1,x_2)$
as $t$ goes to infinity? A more difficult but more natural and
important question is: how to define the hyperbolic gradient flow
for a family of close surfaces? If so, what is the asymptotic
behaviour of a close surface under ``the hyperbolic gradient flow"?
This problem is related to the theory of multi-dimensional
hyperbolic systems of partial differential equations of first order.
}

{\bf 3. The evolution of hypersurfaces in $\mathbb{R}^n\;(n\ge 4)$.}
{\it Investigate the hyperbolic gradient flow in multi-dimensional
Euclidean space $\mathbb{R}^{n+1}\;(n\ge 4)$. In particular, how can
we define a suitable `` hyperbolic gradient flow" to evolve a closed
sub-maifold? if we can, what is the large time behaviour of a close
hypersurface in $\mathbb{R}^n\;(n\ge 4)$ under this kind of
hyperbolic gradient flow. The convex case maybe is easier to study.}

We may also consider variations of the above hyperbolic gradient
flow which can be defined intrinsically on any manifold. For example
we let $(\mathscr{M}, g)$ be a Riemannian manifold, and $X_t\in
\Gamma(\mathscr{M}, TM)$ be a family of tangent vector fileds, the
hyperbolic gradient flow under considered here is given by the
following evolution equation
\begin{equation}\label{5.1}\frac{\partial X_t}{\partial t}+\frac{1}{2}\nabla (\|X_t\|^2)=0,\end{equation}
where, if in local coordinates ${\displaystyle X_t =\sum_{i=1}^n
X^i_t\frac{\partial}{\partial x_i}}$, then $\|X_t\|^2$ is defined by
\begin{equation}\label{5.2}\|X_t\|^2=g_{ij}X_t^iX_t^j\end{equation}
and $\nabla h$ stands for the gradient vector field of a function
$h$ on the manifold, and $g_{ij}= g(\frac{\partial}{\partial
x_i},\frac{\partial}{\partial x_j})$.  By definition, for any given
$h\in C^{\infty}(\mathscr{M}, \mathbb{R})$ and $X\in TM$, we have
\begin{equation}\label{5.3}
g(X, \nabla h)=X(h).\end{equation} The study of this flow will be
very useful to understand the topological and geometrical structure
of the manifold.

Finally, we would like to point out that, perhaps the method in the
present paper is more important than the results obtained here. Our
method may provide a new approach to some conjectures in
differential geometry (see Yau \cite{sy4}).

\vskip 5mm

\noindent{\Large {\bf Acknowledgements.}} This work was completed
while Kong was visiting the Max Planck Institute for Gravitational
Physics (Albert Einstein Institute) during the summer of 2010. Kong
thanks L. Andersson for his invitation and hospitality. This work
was supported in part by the NNSF of China (Grant No. 10971190) and
the Qiu-Shi Chair Professor Fellowship from Zhejiang University,
China.



\begin{thebibliography}{99}

\bibitem{bressan} A. Bressan, Hyperbolic Systems of Conservation Laws:
The One Dimensional Cauchy Problem, Oxford University Press, 2000.

\bibitem{ak} G. Aubert \& P. Kornprobst, Mathematical Problems
in Image Processing,  Springer, 2006.

\bibitem{conway} E. Conway, The formation and decay of shocks for a
 conservation law in several dimensions, {\it Arch. Rat. Mech. Anal.}
 64 (1977), 47-57.

\bibitem{dafermos} C. M. Dafermos, Hyperbolic Conservation Laws in Continuum
Physics, Springer, Berlin Heidelberg, 2005.

\bibitem {dkl0} W.-R. Dai, D.-X. Kong  \& K.-F. Liu, Hyperbolic
geometric flow (I): short-time existence and nonlinear stability,
{\it  Pure and Applied Mathematics Quarterly  (Special Issue: In
honor of Michael Atiyah and Isadore Singer)} 6 (2010), 331-359.

\bibitem {dkl} W.-R. Dai, D.-X. Kong  \& K.-F. Liu, Dissipative hyperbolic
geometric flow, {\it Asian J. Math.} 12 (2008), 345-364.

\bibitem{g} M. Grassin, Global smooth solutions to Euler equations
for a perfect gas, {\it Indiana Univ. Math. J.} 47 (1998),
1397-1432.

\bibitem{h} R. Hamilton, Three-manifolds with positive Ricci curvature, {\it J. Differential Geom.}
17 (1982), 255-306.

\bibitem{he} C.-L. He, D.-X. Kong and K.-F. Liu,
Hyperbolic mean curvature flow, {\it Journal of Differential
Equations} 246 (2009), 373-390.

\bibitem{hi} G. Huisken \& T. Ilmanen, The inverse mean curvature flow and the
Riemannian-Penrose inequality, {\it J. Differential Geom.} 59
(2001), 353-437.

\bibitem {kong} D.-X. Kong, Lectures on Quasilinear Hyperbolic Systems and Applications,
Zhejiang University, Hangzhou, China, 2009.

\bibitem {k} D.-X. Kong,  Hyperbolic geometric flow, {\it the Proceedings of ICCM
2007}, Vol. II, Higher Educationial Press, Beijing, 2007, 95-110.

\bibitem {kl} D.-X. Kong \& K.-F. Liu, Wave character of metrics
and hyperbolic geometric flow, {\it J. Math. Phys.} 48 (2007),
103508.

\bibitem {klw} D.-X. Kong, K.-F. Liu \& Y.-Z. Wang, Life-span of classical solutions to hyperbolic
geometric flow in two space variables with slow decay initial data,
to appear in {\it Communications in Partial Differential Equations}.

\bibitem {kw} D.-X. Kong, K.-F. Liu and Z.-G. Wang,
Hyperbolic mean curvature flow: Evolution of plane curves, {\it Acta
Mathematica Scientia (A special issue dedicated to Professor Wu
Wenjun's 90th birthday)} 29 (2009), 493-514.

\bibitem {klx} D.-X. Kong, K.-F. Liu \& D.-L. Xu, The hyperbolic geometric flow on Riemann surfaces,
{\it Communications in Partial Differential Equations} 34 (2009),
553-580.

\bibitem {kw2} D.-X. Kong and Z.-G. Wang, Formation of singularities in
the motion of plane curves under hyperbolic mean curvature flow,
{\it Journal of Differential Equations} 247 (2009), 1694-1719.

\bibitem {kr}  S. Kruzkov, First-order quasilinear equations with several space
variables, {\it Mathematics of the USSR-Sbornik} 10 (1970), 217-273.

\bibitem{lax} P. D. Lax, Hyperbolic systems of conservation laws II,
{\it Commun. Pure Appl. Math.}  10 (1957), 537-556.

\bibitem {li} T.-T. Li, Global Classical Solutions for Quasilinear
Hyperbolic Systems. RAM: Research in Applied Mathematics, 32.
Masson, Paris; John Wiley \& Sons, Ltd., Chichester, 1994.

\bibitem {p1} G. Perelman, The entropy formula for the Ricci flow and its
geometric applications, {\it arXiv.org}, November 11, 2002.

\bibitem {p2} G. Perelman, Ricci flow with surgery on three-manifolds,
{\it arXiv.org}, March 10, 2003.

\bibitem {p3}  G. Perelman, Finite extinction time for the solutions to the
Ricci flow on certain three-manifolds, {\it arXiv.org}, July 17,
2003.

\bibitem{sy2}  R. Schoen \& S.-T. Yau,  On the proof of the positive mass conjecture in general relativity,
 {\it Comm. Math. Phys.} 65 (1979), 45-76.

\bibitem{sy3}  R. Schoen \& S.-T. Yau, Proof of the positive mass theorem II,
{\it Comm. Math. Phys.} 79 (1981), 231-260.

\bibitem {sy4} R. Schoen \& S.-T. Yau, {\it Lectures on Differential Geometry},
International Press, Cambridge, MA, 1994.

\bibitem{serre} D. Serre, Systems of Conservation Laws 1:
Hyperbolicity, Entropies, Shock Waves, Canmbridge University Press,
Cambridge, 1999.

\end{thebibliography}
\end{document}